\documentclass[a4paper,12pt]{article}

\usepackage{mathrsfs}
\usepackage{amsfonts}
\usepackage{amsmath}
\usepackage{latexsym,bm}
\usepackage{amsthm}
\usepackage{graphicx}
\usepackage[top=1in,bottom=1in,left=1.25in,right=1.25in]{geometry}
\usepackage[dvips]{color}

\newtheorem{theorem}{Theorem}

\newtheorem{lemma}{Lemma}

\newtheorem{remark}{Remark}

\numberwithin{equation}{section} 

\begin{document}

\title{Conditional limit theorems for critical continuous-state branching processes}
\author{Yan-Xia Ren\footnote{Research supported in part by
NNSF of China (Grant No. 10971003 and No. 11128101) and Specialized
Research Fund for the Doctoral Program of Higher Education.}, Ting
Yang and Guo-Huan Zhao}
\date{}
\maketitle

\begin{abstract}
In this paper we study the conditional limit theorems for critical
continuous-state branching processes with branching mechanism
$\psi(\lambda)=\lambda^{1+\alpha}L(1/\lambda)$ where
$\alpha\in [0,1]$ and $L$ is slowly varying at $\infty$. We prove
that if $\alpha\in (0,1]$, there are norming constants $Q_{t}\to 0$
(as $t\uparrow +\infty$) such that for every $x>0$,
$P_{x}\left(Q_{t}X_{t}\in\cdot|X_{t}>0\right)$ converges weakly to a
non-degenerate limit. The converse assertion is also true provided
the regularity of $\psi$ at $0$. We give a conditional limit theorem
for the case $\alpha=0$. The limit theorems we obtain in this paper
allow infinite variance of the branching process.

\end{abstract}

\section{Introduction}
A $[0,+\infty)$-valued strong Markov process $X=\{X_{t}:t\ge 0\}$
with probabilities $\{P_{x}:x>0\}$ is called a (conservative)
continuous-state branching process (CB process) if it has paths that
are right continuous with left limits, and it employs the following
branching property: for any $\lambda\ge 0$ and $x,y>0$,
\begin{equation}
E_{x+y}(e^{-\lambda X_{t}})=E_{x}(e^{-\lambda
X_{t}})E_{y}(e^{-\lambda X_{t}}).\label{branchingproperty}
\end{equation}
It can be characterized by the branching mechanism $\psi$
which is also the Laplace exponent of a L\'{e}vy process with
non-negative jumps. Set $\rho:=\psi'(0+)$, then $E_{x}X_{t}=x
e^{-\rho t}$. We call a CB process \textit{supercritical, critical}
or \textit{subcritical} as $\rho<0,\ =0,$ or $>0$.

Let $\tau:=\inf\{t\ge 0:X_{t}=0\}$ denote the extinction time of
$X_{t}$ and $q(x):=P_{x}(\tau<+\infty)$. When $q(x)<1$ for some (and
then for all) $x>0$, the asymptotic behavior of $X_{t}$ is studied
in \cite{Grey}. It was proved that there are positive constants
$\eta_{t}$ such that $\eta_{t}X_{t}$ converges almost surely to a
non-degenerate random variable as $t\to+\infty$.Note that
$q(x)\equiv1$ if and only if $X$ is subcritical or critical
with $\psi$ satisfying
\begin{equation}
\int_{\theta}^{+\infty}\frac{1}{\psi(\xi)}d\xi<+\infty\label{asump2}
\end{equation}
for some $\theta>0$. In this case, one can study the asymptotic
behavior of $X$ by conditioning it on $\{\tau>t\}$ (see
\cite{Li,Lambert,Pakes2,Pakes3} and the references therein). In the
subcritical case, it was proved that
$P_{x}\left(X_{t}\in\cdot|\tau>t\right)$ converges weakly as
$t\to+\infty$ to the so-called Yaglom distribution. However in the critical case, the
limiting distribution of $X_{t}$ conditioned on non-extinction is
trivial, converging to the Dirac measure at $\infty$. To evaluate
the asymptotic behavior of $X_{t}$ more accurately, we therefore
have to normalize the process appropriately.

Throughout this paper,
we assume $\psi$ satisfies
\begin{equation}
\psi(\lambda)=\lambda^{1+\alpha}L(1/\lambda)\quad \forall \lambda\ge
0\label{psi}
\end{equation}
where $\alpha\in [0,1]$ and $L$ is slowly varying at infinity. Our assumption
on $\psi$ does not require the finiteness of $E_{x}X_{t}^{2}$.

It is well known that a CB process can be viewed as the analogue of
Galton-Watson branching process in continuous time and continuous
state space. So it is necessary for us to take a look at the
asymptotic behavior of critical G-W branching processes. Let $f(s)$
denote the probability generating function of the offspring law of
the critical G-W process $Z_{n}$. Let $\bar{F}(n)=P_{1}(Z_{n}>0)$.
Slack \cite{Slack1,Slack2} proved that $P_{1}( \bar{F}(n)Z_{n}\le
y|Z_{n}>0)$ converges weakly to a non-degenerate limit if and only
if
\begin{equation}
f(s)=s+(1-s)^{1+\alpha}L\left(\frac{1}{1-s}\right)\label{fs}
\end{equation}
for some $\alpha\in (0,1]$ and $L$ slowly varying at $+\infty$.
Later Nagaev \textit{et.al.}\cite{Nagaev} proved a conditional limit
theorem for $f(s)$ satisfying \eqref{fs} with $\alpha=0$. Recently,
Pakes \cite{Pakes1} generalized the above results to continuous
time Markov branching process. The proofs given in \cite{Pakes1},
based on Karamata's theory for regular varying functions,
are much easier.
However, for discrete-state branching process, there leaves open the
question of whether \eqref{fs} is implied by the more general
conditional convergence of $P_{1}(b_{n}Z_{n}\le y|Z_{n}>0)$ for some
positive sequence $\{b_{n}\}$ with $b_{n}\to 0$.

This paper is structured as follows: In Section 2, we collect some
basic facts about regularly varying functions and CB processes.
Section 3 is devoted to the conditional limit theorems for $\psi$
with $\alpha\in (0,1]$. We prove that there exists positive norming
constants $Q_{t}\to 0$ such that $P_{x}(Q_{t}X_{t}\in \cdot|\tau>t)$
converges weakly to a non-degenerate limit. An admissible norming is
$Q_{t}=P_{1}(\tau>t)$. This is analogous to the result we mentioned
in the above paragraph for discrete-state branching processes. Later
we prove that the converse assertion is also true provided some
regularity of $\psi$ at $0$ (or equivalently, provided some
regularity of the L\'{e}vy measure of $\psi$ at infinity). In Section 4,
we give a conditional limit theorem for the case
$\alpha=0$. Its discrete state analogue is proved independently in
\cite{Nagaev} and \cite{Pakes1}. The last section provides some
concrete examples which satisfy the assumptions in Section 3 or Section 4.
The branching mechanisms in these examples are well
known and taken from \cite{Ren}.

\section{Preliminary}

In the rest of this paper, we shall use the notation $f(x)\sim g(x)$
for functions $f$ and $g$ to mean that $f(x)/g(x)\to 1$ as $x\to +\infty$ or $0$.
Let $x\wedge y:=\min\{x,y\}$.

Suppose $X$ is a CB process with branching mechanism $\psi$. Generally $\psi$ is specified by the L\'{e}vy-Khintschine
formula
$$\psi(\lambda)=a\lambda+b\lambda^{2}+\int_{(0,+\infty)}(e^{-\lambda x}-1+\lambda x)\Lambda(dx),\quad \lambda\ge 0,$$
where $a\in (-\infty,+\infty)$, $b\ge 0$ and $\Lambda$ is a
non-negative measure on $(0,+\infty)$ satisfying
$\int_{(0,+\infty)}(x^{2}\wedge x)\Lambda(dx)<+\infty$. $\Lambda$ is
called the L\'{e}vy measure of $\psi$. Obviously, $\psi$ is convex
and infinitely differentiable on $(0,+\infty)$. Since we aim at
conditioning critical CB process on non-extinction, we assume that
$\psi$ satisfies \eqref{asump2} with $\psi'(0+)=0$. Under this
assumption, $\psi$ is a strictly convex function on $[0,+\infty)$,
$\psi(+\infty)=+\infty$, and $\psi(\lambda)=0$ if and only if
$\lambda=0$. This assumption also implies that
$P_{x}(\tau<+\infty)=1$ for every $x>0$.

For $x>0$ and $\lambda,t\ge 0$, let $E_{x}(e^{-\lambda X_{t}})=e^{-x
u_{t}(\lambda)}$. Then $u_{t}(\lambda)$ is the unique positive
solution to the backward equation
\begin{equation}
\frac{\partial}{\partial
t}u_{t}(\lambda)=-\psi(u_{t}(\lambda)),\quad
u_{0}(\lambda)=\lambda.\label{backward}
\end{equation}
From $\eqref{backward}$ and the semi-group property
$u_{t}(u_{s}(\lambda))=u_{t+s}(\lambda)$, we also get the forward
equation
\begin{equation}
\frac{\partial}{\partial
t}u_{t}(\lambda)=-\psi(\lambda)\frac{\partial}{\partial
\lambda}u_{t}(\lambda),\quad u_{0}(\lambda)=\lambda.\label{forward}
\end{equation}
 Note that our moment condition on $\Lambda$
implies that $E_{x}X_{t}=x e^{-\rho t}<+\infty$ for all $x>0$ and
$t\ge 0$.

Next define
$$\phi(z):=\int_{z}^{+\infty}\frac{1}{\psi(\xi)}\,d\xi,\quad \forall z>0.$$
The mapping $\phi:(0,+\infty)\to (0,+\infty)$ is bijective with
$\phi(0)=+\infty$ and $\phi(+\infty)=0$. We use $\varphi$ to denote
the inverse function of $\phi$. From \eqref{backward}, we have
$$\int_{u_{t}(\lambda)}^{\lambda}\frac{1}{\psi(\xi)}\,d\xi=t,\quad \lambda,t\ge 0.$$
Hence
\begin{equation}
u_{t}(\lambda)=\varphi(t+\phi(\lambda)),\quad \lambda,t\ge
0.\label{eq1}
\end{equation}
Since $\phi(+\infty)=0$, we have $u_{t}(+\infty)=\varphi(t)$, and
for any $x>0$ and $t\ge 0$,
\begin{equation}
P_{x}(\tau>t)=P_{x}(X_{t}>0)=1-\lim_{\lambda\to+\infty}e^{-x
u_{t}(\lambda)}=1-e^{-x\varphi(t)}.\label{extinprob}
\end{equation}
Let $\bar{F}(t):=P_{1}(\tau>t)$. Obviously, we have $\bar{F}(t)\sim
\varphi(t)$ as $t\uparrow +\infty$.

Results about regular varying functions
will be used a lot in the remaining
paper, so we collect some basic facts here. A positive measurable
function $L$ is said to be slowly varying at $\infty$ if it is
defined on $(0,+\infty)$ and $\lim_{x\to+\infty}L(\lambda x)/L(x)=1$
for all $\lambda>0$. This convergence holds uniformly with respect
to $\lambda$ on every compact subset of $(0,+\infty)$. Let
$\mathcal{S}$ denote the set of all slowly varying functions at
$\infty$. If $L\in \mathcal{S}$, then for any $\delta>0$,
$\lim_{x\to+\infty}x^{\delta}L(x)=+\infty$, and
$\lim_{x\to+\infty}x^{-\delta}L(x)=0$.

If a positive function $f$ defined on $(0,+\infty)$ satisfies that
$f(\lambda x)/f(x)\to\lambda^{p}$ as $x\to+\infty$ (resp. $0$) for
any $\lambda>0$, then $f$ is called regularly varying at $\infty$
(resp. $0$) with index $p\in(-\infty,+\infty)$, denoted by $f\in
\mathcal{R}_{p}(\infty)$ (resp. $f\in \mathcal{R}_{p}(0)$).
Obviously, $f(x)\in \mathcal{R}_{p}(0)$ is equivalent to $f(1/x)\in
\mathcal{R}_{-p}(\infty)$. If $f\in\mathcal{R}_{p}(\infty)$ (resp.
$f\in\mathcal{R}_{p}(0)$), it can be represented by $f(x)=x^{p}L(x)$
(resp. $f(x)=x^{p}L(1/x)$) for some $L\in \mathcal{S}$.

\section{The case $0<\alpha\le 1$}

The following technical lemma follows from Theorem 1.5.2 and Theorem 1.5.12 in \cite{Bingham}.
We omit the details here.

\begin{lemma}\label{lem1}
\begin{enumerate}
\item[(1)] If $p\in(-\infty,+\infty)$, $f\in \mathcal{R}_{p}(\infty)$ (resp. $\mathcal{R}_{p}(0)$),
$T_{1}(t),T_{2}(t)\to +\infty$ (resp. $0$) and $T_{1}(t)\sim
T_{2}(t)$ as $t\uparrow +\infty$, then $f(T_{1}(t))\sim
f(T_{2}(t))$.

\item[(2)] Suppose $f\in\mathcal{R}_{p}(\infty)$, $T_{1}(t),T_{2}(t)\to
+\infty$ as $t\to +\infty $, and $f(T_{1}(t))/f(T_{2}(t))\sim c\in
(0,+\infty)$. If $p>0$, then $T_{1}(t)/T_{2}(t)\sim c^{1/p}$;
otherwise if $p<0$ and $f$ has inverse function $f^{-1}$, then
$f^{-1}\in \mathcal{R}_{1/p}(0)$ and $T_{1}(t)/T_{2}(t)\sim
c^{1/p}$.
\end{enumerate}
\end{lemma}

\begin{theorem}\label{them1}
If \eqref{psi} holds with $0<\alpha\le 1$, then for all $x>0$ and
$y\ge 0$,
\begin{equation}
\lim_{t\to+\infty}P_{x}\left(\bar{F}(t)X_{t}\le y |
\tau>t\right)=H_{\alpha}(y),\label{1.eq1}
\end{equation}
where $H_{\alpha}(y)$ is a probability distribution function, and
its Laplace transform is given by
\begin{equation}
h_{\alpha}(\theta)=\int_{[0,+\infty)}e^{-\theta
y}dH_{\alpha}(y)=1-(1+\theta^{-\alpha})^{-1/\alpha}.\label{1.eq2}
\end{equation}
Moreover, $\bar{F}(t)$ is regularly varying at $+\infty$ with index
$-1/\alpha$, and consequently, for any $\delta>0$,
$$\lim_{t\to+\infty}t^{\frac{1}{\alpha}+\delta}\bar{F}(t)=+\infty,\quad \lim_{t\to+\infty}t^{\frac{1}{\alpha}-\delta}\bar{F}(t)= 0.$$
\end{theorem}

\proof For any $z>0$, set
$g(z):=\phi(1/z)=\int_{0}^{z}\xi^{\alpha-1}/L(\xi)\,d\xi.$
Then by Karamata's theorem (see, for example \cite[Theorem
1.5.11]{Bingham}), we have $g\in \mathcal{R}_{\alpha}(\infty)$, more
specifically, $g(z)\sim \alpha^{-1}z^{\alpha}L(z)^{-1}$ as $z\to
+\infty$. Consequently, we get $\phi\in\mathcal{R}_{-\alpha}(0)$,
$\phi(z)\sim \alpha^{-1}z^{-\alpha}L(1/z)^{-1}$ as $z\downarrow 0$,
and $\varphi \in \mathcal{R}_{-1/\alpha}(\infty)$.

Since $1-e^{-u}\sim u$ as $u\downarrow 0$, we have for any
$x,\theta>0$,
\begin{eqnarray}
\lim_{t\to+\infty}E_{x}\left(e^{-\theta\bar{F}(t)X_{t}}|\tau>t\right)&=&1-\lim_{t\to+\infty}\frac{1-e^{-x\varphi(t+\phi(\theta\bar{F}(t)))}}{1-e^{-x\varphi(t)}}\nonumber\\
&=&1-\lim_{t\to
+\infty}\frac{\varphi(t+\phi(\theta\bar{F}(t)))}{\varphi(t)}.\label{them1.1}
\end{eqnarray}
It follows from Lemma \ref{lem1} and the fact that $\bar{F}(t)\sim \varphi(t)$ as
$t\uparrow +\infty$, we have
$$\phi(\theta\bar{F}(t))\sim \phi(\theta\varphi(t))\sim
\theta^{-\alpha}\phi(\varphi(t))=\theta^{-\alpha}t.$$ Hence we have
$\varphi(t+\phi(\theta\bar{F}(t)))\sim\varphi((1+\theta^{-\alpha})t)$.
By \eqref{them1.1} and the regularity of $\varphi$ at $\infty$, we
get
\begin{equation}
\lim_{t\to+\infty}E_{x}\left(e^{-\theta\bar{F}(t)X_{t}}|\tau>t\right)=1-\lim_{t\to+\infty}\frac{\varphi((1+\theta^{-\alpha})t)}{\varphi(t)}=1-(1+\theta^{-\alpha})^{-1/\alpha}.\label{them1.2}
\end{equation}
The assertion follows from the continuity theory for Laplace transforms (see, for example, \cite[Section 6.6 ]{Chung}). \qed

\begin{remark} The stationary-excess operation on $H_{\alpha}(y)$ is defined by
$\widetilde{H}_{\alpha}(y):=\int_{(0,y]}\bar{H}_{\alpha}(x)dx/\int_{(0,+\infty)}\bar{H}_{\alpha}(x)dx$,
where $\bar{H}_{\alpha}(y)=1-H_{\alpha}(y)$.
$\widetilde{H}_{\alpha}(y)$ is also a probability distribution
function, and a simple calculation shows that its Laplace transform
is $(1+\theta^{-\alpha})^{-1/\alpha}$. $\widetilde{H}_{\alpha}(y)$ is
often called a generalized positive Linnik law. When $\alpha=1$, it
gives the well-known standard exponential law. For more information
on Linnik Law, we refer readers to \cite[Section 4]{Pakes1} and references therein.
\end{remark}

\bigskip

The remainder of this section is devoted to the converse assertions
to Theorem \ref{them1}. Suppose that $X_{t}$ is a critical CB process. If there exist
$x>0$ and positive constants $Q_{t}\to 0$ (as $t\uparrow +\infty$)
such that $P_{x}\left(Q_{t}X_{t}\in\cdot|\tau>t\right)$ converges
weakly to a non-degenerate limit, then
$\liminf_{t\to+\infty}Q_{t}/\bar{F}(t)>0$. In fact, by
Fatou's lemma
\begin{eqnarray*}0&<&\liminf_{t\to+\infty}\int_{0}^{+\infty}P_{x}\left(Q_{t}X_{t}>y|\,\tau>t\right)dy\\
&=&\liminf_{t\to+\infty}E_{x}\left(Q_{t}X_{t}\,|\,\tau>t\right)\\
&=&\liminf_{t\to+\infty}Q_{t}/\bar{F}(t).\end{eqnarray*}

\begin{lemma}\label{lem2}
Suppose $\psi$ is the branching mechanism of a non-trivial critical CB process. If $\psi$ is regularly varying at $0$, then
$\psi\in \mathcal{R}_{1+\alpha}(0)$ with $\alpha\in [0,1]$.
\end{lemma}

\proof Suppose $\psi(\lambda)=\lambda^{p}L(1/\lambda)$ for some
$p\in(-\infty,+\infty)$ and $L\in\mathcal{S}$. Since
$$0=\psi'(0+)=\lim_{\lambda\downarrow
0}\frac{\psi(\lambda)}{\lambda}=\lim_{\lambda\downarrow
0}\lambda^{p-1}L(1/\lambda),$$ we have $p\ge 1$. If $p>2$, then
\begin{equation}
\psi''(0+)=\lim_{\lambda\downarrow
0}\frac{2\psi(\lambda)}{\lambda^{2}}=\lim_{\lambda\downarrow
0}2\lambda^{p-2}L(1/\lambda)=0.\label{lem2.1}
\end{equation}
Recall that $\psi''(\lambda)=2b+\int_{0}^{+\infty}x^{2}e^{-\lambda
x}\Lambda(dx)$ for some $b\ge 0$ and $\int_{(0,+\infty)}(x\wedge
x^{2})\Lambda(dx)<+\infty$. So \eqref{lem2.1} implies that $b=0$
and $\Lambda(dx)\equiv 0$, in which case $\psi$ is trivial. Hence
$p\le 2$. We set $\alpha=p-1$, thus proving the conclusion.\qed

\begin{theorem}\label{them2}
Suppose $X_{t}$ is a critical CB process with branching mechanism
$\psi$. If for some $x>0$, $P_{x}\left(\bar{F}(t)X_{t}\le y |
\tau>t\right)$ converges weakly to a non-degenerate distribution
function $H(y)$, then \eqref{psi} holds with $\alpha\in (0,1]$.
\end{theorem}
\proof Let $H(y,t):=P_{x}\left(\bar{F}(t)X_{t}\le y |
\tau>t\right)$. Under the assumption, we have
\begin{equation}\label{weak-limit}\lim_{t\to+\infty}\int_{[0,+\infty)}g(y)dH(y,t)=\int_{[0,+\infty)}g(y)dH(y)\end{equation}
for any continuous function $g$ defined on $[0,+\infty)$ such that
$\lim_{y\to+\infty}g(y)=0$. Suppose $\theta>0$.
Using \eqref{weak-limit} with $g(y)=e^{-\theta y}$ we get
\begin{eqnarray}
h(\theta):=\int_{[0,+\infty)}e^{-\theta
y}dH(y)&=&\lim_{t\to+\infty}\int_{[0,+\infty)}e^{-\theta
y}dH(y,t)\nonumber\\
&=&\lim_{t\to+\infty}E_{x}\left(e^{-\theta
\bar{F}(t)X_{t}}|\tau>t\right)\nonumber\\
&=&1-\lim_{t\to+\infty}\frac{1-\exp\{-x
u_{t}(\theta\bar{F}(t))\}}{1-\exp\{-x\varphi(t)\}}\nonumber\\
&=&1-\lim_{t\to+\infty}\frac{u_{t}(\theta\bar{F}(t))}{\varphi(t)}.
\end{eqnarray}
So as $t\uparrow +\infty$
\begin{equation}
u_{t}(\theta\bar{F}(t))\sim\bar{h}(\theta)\varphi(t)\sim
\bar{h}(\theta)\bar{F}(t),\label{them2.1}
\end{equation}
where $\bar{h}(\theta)=1-h(\theta)$. On the other hand,using \eqref{weak-limit} with $g(y)=y\,e^{-\theta y}$, we obtain
\begin{eqnarray}
\bar{h}'(\theta)=\int_{[0,+\infty)}y e^{-\theta
y}dH(y)&=&\lim_{t\to+\infty}\int_{[0,+\infty)}y
e^{-\theta y}dH(y,t)\nonumber\\
&=&\lim_{t\to
+\infty}E_{x}\left(\bar{F}(t)X_{t}e^{-\theta\bar{F}(t)X_{t}}|\tau>t\right)\nonumber\\
&=&\lim_{t\to+\infty}\frac{\bar{F}(t)E_{x}(X_{t}e^{-\theta\bar{F}(t)X_{t}})}{1-e^{-x\varphi(t)}}.\label{them2.2}
\end{eqnarray}

From \eqref{backward} and \eqref{forward}, we have
$$\frac{\partial}{\partial \lambda}u_{t}(\lambda)=\frac{\psi(u_{t}(\lambda))}{\psi(\lambda)},\quad\forall \lambda> 0.$$
Thus
\begin{equation}
E_{x}(X_{t}e^{-\lambda X_{t}})=-\frac{\partial}{\partial
\lambda}e^{-x u_{t}(\lambda)}=xe^{-x
u_{t}(\lambda)}\frac{\psi(u_{t}(\lambda))}{\psi(\lambda)}.\label{them2.3}
\end{equation}
It follows from \eqref{them2.1}, \eqref{them2.2} and\eqref{them2.3} that
\begin{eqnarray}
\bar{h}'(\theta)&=&\lim_{t\to+\infty}\frac{x\bar{F}(t)}{1-e^{-x\varphi(t)}}\,e^{-x
u_{t}(\theta\bar{F}(t))}\frac{\psi(u_{t}(\theta\bar{F}(t)))}{\psi(\theta\bar{F}(t))}\nonumber\\
&=&\lim_{t\to+\infty}\frac{\psi(u_{t}(\theta\bar{F}(t)))}{\psi(\theta\bar{F}(t))}\nonumber\\
&=&\lim_{t\to+\infty}\frac{\psi(\bar{h}(\theta)\bar{F}(t))}{\psi(\theta\bar{F}(t))}.\label{them2.4}
\end{eqnarray}
The last equality follows from a standard argument using the
continuity and monotonicity of $\psi$. Let
$\lambda(\theta):=\bar{h}(\theta)/\theta=\int_{0}^{+\infty}e^{-\theta
y}\bar{H}(y)dy$ where $\bar{H}(y)=1-H(y)$. $\lambda(\theta)$ is
decreasing on $(0,+\infty)$. Since $\bar{F}(t)$ decreases
continuously to $0$ as $t\uparrow +\infty$ and $\psi$ is monotone on
$(0,+\infty)$, \eqref{them2.4} implies that
\begin{equation}
\lim_{s\downarrow
0}\frac{\psi(\lambda(\theta)s)}{\psi(s)}=\xi(\lambda(\theta)),\quad\forall
\theta>0,\label{them2.5}
\end{equation}
for some function $\xi$ such that
$\xi(\lambda(\theta))=\bar{h}'(\theta)$. From the continuity and
monotonicity of $\lambda(\theta)$, we have for any $\lambda\in
(0,\lambda(0+))$,
\begin{equation}
\lim_{s\downarrow 0}\frac{\psi(\lambda
s)}{\psi(s)}=\xi(\lambda).\label{them2.6}
\end{equation}
Characterization theorem (see \cite[Theorem 1.4.1]{Bingham} ) says that
\eqref{them2.6} holds for all $\lambda>0$, and there exists $p\in
(-\infty,+\infty)$ such that $\xi(\lambda)\equiv\lambda^{p}$, \textit{i.e.}
$\psi$ is regularly varying at $0$ with index $p$. Let $\alpha=p-1$,
then $\alpha\in [0,1]$ by Lemma \ref{lem2}. If $\alpha=0$, we have
$$\frac{\bar{h}(\theta)}{\theta}=\lambda(\theta)=\xi(\lambda(\theta))=\bar{h}'(\theta).$$
This has the solution $h(\theta)=1-c\theta$ for some constant $c$.
This is the Laplace transform of a distribution function if and only
if $c=0$, in which case $H(y)\equiv 1$ is the distribution function
of Dirac measure at $0$. Therefore $\alpha>0$.\qed

\bigskip

Suppose $\mu$ is a positive measure supported on $(0,+\infty)$. We
say $\mu$ is regularly varying at $+\infty$ if $u(x):=\mu((0,x])$ is
regularly varying at $+\infty$. The following theorem tells us that
\eqref{psi} with $\alpha\in (0,1]$ is implied by the more general
limit $P_{x}\left(Q_{t}X_{t}\le y|\tau>t\right)\to H(y)$ where
$Q_{t}$ are positive constants such that $Q_{t}\to 0$.

\begin{theorem}\label{them3}
Let $\psi$ be the branching mechanism of a non-trivial
critical CB process with L\'{e}vy measure $\Lambda$. Suppose
$x^{2}\Lambda(dx)$ is regularly varying at $+\infty$. If there exist
$x>0$ and positive constants $Q_{t}\to 0$ (as $t\uparrow +\infty$)
such that $P_{x}\left(Q_{t}X_{t}\le y|\tau>t\right)$ converges
weakly to a non-degenerate limit $H(y)$, then \eqref{psi} holds with
$\alpha\in (0,1]$. In this case, $Q_{t}/\bar{F}(t)\sim c\in
(0,+\infty)$, and the Laplace transform of $H(y)$ is given by
$$h(\theta)=\int_{[0,+\infty)}e^{-\theta y}dH(y)=1-(1+c^{-\alpha}\theta^{-\alpha})^{-1/\alpha}.$$
\end{theorem}

To proof Theorem \ref{them3}, we need the following lemma.
\begin{lemma}\label{prop1}
Suppose $\psi$ is the branching mechanism of a non-trivial critical CB process. Then $\psi$ is regularly varying at $0$ if and
only if $x^{2}\Lambda(dx)$ is regularly varying at $+\infty$.
\end{lemma}
\proof We may and do assume that
$$\psi(\lambda)=b\lambda^{2}+\int_{(0,+\infty)}(e^{-\lambda x}-1+\lambda x)\Lambda(dx)$$
where $b\ge 0$ and $\int_{(0,+\infty)}(x\wedge
x^{2})\Lambda(dx)<+\infty$. Let $U(z):=\int_{(0,z]}x^{2}\Lambda(dx)$
and $\hat{U}(\theta):=\int_{(0,+\infty)}e^{-\theta x}dU(x)$. If
$\psi''(0+)<+\infty$, then $\psi\in\mathcal{R}_{2}(0)$ and
$\int_{[1,+\infty)}x^{2}\Lambda(dx)<+\infty$. Obviously
$\lim_{z\to+\infty}U(z)=\int_{(0,+\infty)}x^{2}\Lambda(dx)<+\infty$,
which implies that $x^{2}\Lambda(dx)$ is slowly varying at
$+\infty$.

Now we suppose $\psi''(0+)=+\infty$, in which case
$\int_{[1,+\infty)}x^{2}\Lambda(dx)=+\infty$. If $\psi$ is regularly
varying at $0$ with index $p\in [1,2]$, then for any $A>0$, using
L'Hospital rule, we have
\begin{eqnarray}
A^{p}&=&\lim_{\lambda\to
0+}\frac{\psi(A\lambda)}{\psi(\lambda)}=\lim_{\lambda\to
0+}A^{2}\,\frac{\psi''(A\lambda)}{\psi''(\lambda)}\nonumber\\
&=&\lim_{\lambda\to
0+}A^{2}\,\frac{2b+\hat{U}(A\lambda)}{2b+\hat{U}(\lambda)}=\lim_{\lambda\to
0+}A^{2}\,\frac{\hat{U}(A\lambda)}{\hat{U}(\lambda)}.\label{prop1.1}
\end{eqnarray}
The last equality is because $\lim_{\theta\to
0+}\hat{U}(\theta)=\lim_{\theta\to 0+}\int_{(0,+\infty)}e^{-\theta
x}x^{2}\Lambda(dx)=+\infty$. Thus $\hat{U}$ is regularly varying at
$0$ with index $p-2\in [-1,0]$. By Tauberian theorem (see, for
example \cite[Theorem 1.7.1]{Bingham}), $x^{2}\Lambda(dx)$ is
regularly varying at $+\infty$ with index $2-p\in [0,1]$. The
converse assertion is clear through the equalities in
\eqref{prop1.1}.\qed

\bigskip

\textit{Proof of Theorem \ref{them3}.} The proof is similar to that
of Theorem \ref{them2}. We provide details here for the reader's
convenience. Let $H(y,t):=P_{x}\left(Q_{t}X_{t}\le y|\tau>t\right)$,
$h(\theta):=\int_{[0,+\infty)}e^{-\theta y}dH(y,t)$ and
$\bar{h}(\theta):=1-h(\theta)$. Similarly we can get the analogues
to \eqref{them2.1} and \eqref{them2.4}:
\begin{equation}
u_{t}(\theta Q_{t})\sim \bar{h}(\theta)\bar{F}(t)\ \mbox{ as
}t\to+\infty,\label{them3.1}
\end{equation}
and
\begin{equation}
\lim_{t\to+\infty}\frac{Q_{t}}{\bar{F}(t)}\,\frac{\psi(u_{t}(\theta
Q_{t}))}{\psi(\theta Q_{t})}=\bar{h}'(\theta).\label{them3.2}
\end{equation}
It follows from Lemma \ref{prop1} that $\psi$ is regularly varying
at $0$. Using Lemma \ref{lem1}, \eqref{them3.1} and \eqref{them3.2},
we have
\begin{equation}
\lim_{t\to+\infty}\frac{Q_{t}}{\bar{F}(t)}\,\frac{\psi(\bar{h}(\theta)\bar{F}(t))}{\psi(\theta
Q_{t})}=\bar{h}'(\theta).\label{them3.3}
\end{equation}

In view of Lemma \ref{lem2}, we may and do assume
$\psi\in\mathcal{R}_{1+\alpha}(0)$ with $\alpha\in [0,1]$. We first
consider the case $\alpha>0$. Put $g(z):=(z\psi(1/z))^{-1}$, $z>0$.
Then $g\in \mathcal{R}_{\alpha}(+\infty)$.
\eqref{them3.3} implies that
\begin{equation}
 \lim_{t\to+\infty}\frac{g(1/\theta Q_{t})}{g(1/\bar{h}(\theta)\bar{F}(t))}
=\lim_{t\to+\infty}\frac{\psi(\bar{h}(\theta)\bar{F}(t))}{\psi(\theta
Q_{t})}\,\frac{\theta
Q_{t}}{\bar{h}(\theta)\bar{F}(t)}=\frac{\theta}{\bar{h}(\theta)}\,\bar{h}'(\theta),\quad\forall
\theta>0.\label{them3.4}
\end{equation}
By Lemma \ref{lem1}, we have for all $\theta>0$,
$$\frac{\theta Q_{t}}{\bar{h}(\theta)\bar{F}(t)}\sim \left(\frac{\theta}{\bar{h}(\theta)}\,\bar{h}'(\theta)\right)^{-1/\alpha},\quad \mbox{as } t\uparrow +\infty,$$
or equivalently,
$$\frac{Q_{t}}{\bar{F}(t)}\sim \left(\frac{\theta}{\bar{h}(\theta)}\right)^{-1/\alpha-1}\bar{h}'(\theta)^{-1/\alpha},\quad \mbox{as } t\uparrow +\infty.$$
Hence we have $Q_{t}/\bar{F}(t)\sim c$ for some constant $c\in
(0,+\infty)$, and
$$\left(\frac{\theta}{\bar{h}(\theta)}\right)^{-1/\alpha-1}\bar{h}'(\theta)^{-1/\alpha}
\equiv c,\quad \theta\in(0,\infty).$$
In view of the initial condition $\bar{h}(0)=1$, the above equation has the unique solution
$h(\theta)=1-(1+c^{-\alpha}\theta^{-\alpha})^{-1/\alpha}$.

Otherwise if $\alpha=0$, we assume $\psi(\lambda)=\lambda
l(\lambda)$ where $l$ is slowing varying at $0$. From
\eqref{them3.3}, we get
$$\lim_{t\to +\infty}\frac{l(\bar{F}(t))}{l(Q_{t})}=\frac{\theta}{\bar{h}(\theta)}\,\bar{h}'(\theta),\quad\forall \theta>0.$$
Thus there exists a constant $c_{1}$ independent of $\theta$ such
that
$$\frac{\theta}{\bar{h}(\theta)}\,\bar{h}'(\theta)
\equiv c_{1},\quad \theta\in(0, \infty).$$
This has the solution
$h(\theta)=1-c_{2}\theta^{c_{1}}$ for some constant $c_{2}$.
$h(\theta)$ is the Laplace transform of a distribution function only if $c_{2}=0$,
in which case $H(y)\equiv 1, y\in [0, \infty)$ is the distribution
function of the Dirac measure at $0$. This contradicts our assumption that $H$ is the distribution function of a non-degenerate random variable.
Hence $\alpha>0$. We complete the proof.\qed

\begin{remark}\label{remark2}
Through the above proof we see that for $\psi$ satisfying
\eqref{psi} with $\alpha=0$, the limit distribution of
$P_{x}\left(Q_{t}X_{t}\in\cdot\,|\,\tau>t\right)$, if exists, must
be the Dirac measure at $0$.
\end{remark}

\section{The case $\alpha=0$}

In this section, we stay in the regime $\alpha=0$. Suppose
$\psi(\lambda)=\lambda L(1/\lambda)$ satisfies our assumption
\eqref{asump2} and $\psi'(0+)=0$. From Remark \ref{remark2} we know
that for $\alpha=0$, any possible positive sequence $Q_{t}\to 0$
overnormalizes $X_{t}$. So we need to find an alternative way to
normalize $X_{t}$. \cite{Pakes1} considers the analogous conditional
limit theorem for critical Markov branching processes with the
offspring generating function $f(s)=s+(1-s)L(1/(1-s))$ where $L\in
\mathcal{S}$. The proof in \cite{Pakes1} can be adapted here to get
the convergence result for a CB process.

Set
$$V(x):=\phi(1/x)=\int_{1/x}^{+\infty}\frac{1}{\psi(\xi)}\,d\xi=\int_{0}^{x}\frac{1}{\xi L(\xi)}\,d\xi,\quad x>0.$$
Obviously, $V$ is differentiable, strictly increasing on
$(0,+\infty)$, $V'(x)=x^{-1}L(x)^{-1}$, $V(0)=0$ and
$V(+\infty)=\int_{0}^{+\infty}1/\psi(\xi)d\xi=+\infty$. By
Karamata's theorem, we have $V\in \mathcal{S}$, and $V(x)L(x)\to
+\infty$ as $x\to+\infty$.

Let $R$ denote the inverse function of $V$. It is easy to see that
$R(x)=1/\varphi(x)$, $R$ is continuous, strictly increasing on
$(0,+\infty)$ with $R(+\infty)=+\infty$ and $R(0)=0$. By \cite[Theorem 2.4.7]{Bingham},
$R$ belongs to the class of Karamata
rapidly varying functions denoted by $KR_{\infty}$. We refer readers
to \cite[Section 2.4]{Bingham} for more information about
$KR_{\infty}$. Since $y=V(R(y))$, we have
$$1=V'(R(y))R'(y)=\frac{R'(y)}{R(y)L(R(y))},\quad\forall y>0,$$
or equivalently
$$\frac{R'(y)}{R(y)}=L(R(y)),\quad\forall y>0.$$
Thus there exist $c,A>0$ such that
\begin{equation}
R(y)=c\exp\left\{\int_{A}^{y}L(R(z))dz\right\}\label{R(y)},\quad y\in [A,+\infty).
\end{equation}

\begin{lemma}[\cite{Pakes1} Lemma 5.2]\label{lem3}
As $t\uparrow +\infty$, $I(y,t):=\int_{t}^{t+y/L(R(t))}L(R(z))dz\to
y$, and this convergence holds locally uniformly with respect to
$y\in (-\infty,+\infty)$.
\end{lemma}

\begin{theorem}\label{them4}
If \eqref{psi} holds with $\alpha=0$, then
\begin{equation}
V(\bar{F}(t)^{-1})\sim t,\mbox{\quad as }t\uparrow
+\infty,\label{4.eq1}
\end{equation}
and
\begin{equation}
\lim_{t\to +\infty}P_{x}\left(L(\bar{F}(t)^{-1})V(X_{t})\le y
|\tau>t\right)=1-e^{-y}\label{4.eq2}
\end{equation}
for any $x>0$ and $y\ge 0$.
\end{theorem}

\proof \eqref{4.eq1} follows from the fact that $V(\bar{F}(t)^{-1})\sim V(R(t))=t$ as $t\uparrow +\infty$.
Henceforth we only need to prove
\eqref{4.eq2}. By the monotonicity of $V$, we have
\begin{equation}
P_{x}\left(L(\bar{F}(t)^{-1})V(X_{t})\le y
|\tau>t\right)=P_{x}\left(X_{t}\le
R\left(y/L(\bar{F}(t)^{-1})\right) |\tau>t\right).\label{them4.0}
\end{equation}
For any $\theta>0$, using the argument of \eqref{them1.1},
we have
\begin{eqnarray}
&&\lim_{t\to+\infty}P_{x}\left(\exp\left\{-\theta\,\frac{X_{t}}{R\left(y/L(\bar{F}(t)^{-1})\right)}\right\}|\tau>t\right)\nonumber\\
&=&1-\lim_{t\to+\infty}\frac{\varphi\left(t+\phi\left(\theta/R\left(y/L(\bar{F}(t)^{-1})\right)\right)\right)}{\varphi(t)}\nonumber\\
&=&1-\lim_{t\to+\infty}\frac{R(t)}{R\left(t+\phi\left(\theta/R\left(y/L(\bar{F}(t)^{-1})\right)\right)\right)},\label{them4.1}
\end{eqnarray}
where in the last equality we used the fact that $R(t)=1/\varphi(t), t>0$.

Since $V\in \mathcal{S}$ and $\bar{F}(t)\sim \varphi(t)=R(t)^{-1}$
as $t\uparrow +\infty$, we get
\begin{eqnarray}
\phi\left(\theta/R\left(y/L(\bar{F}(t)^{-1})\right)\right)&=&V\left(\frac{1}{\theta}R\left(y/L(\bar{F}(t)^{-1})\right)\right)\nonumber\\
&\sim&
V\left(R\left(y/L(\bar{F}(t)^{-1})\right)\right)\nonumber\\
&=&\frac{y}{L(\bar{F}(t)^{-1})}\nonumber\\
&\sim& \frac{y}{L(R(t))}.\label{them4.2}
\end{eqnarray}
 Thus by \eqref{R(y)}, \eqref{them4.2} and Lemma \ref{lem3}, we have
\begin{eqnarray}
&&\lim_{t\to+\infty}\frac{R(t)}{R\left(t+\phi\left(\theta/R\left(y/L(\bar{F}(t)^{-1})\right)\right)\right)}\nonumber\\
&=&\lim_{t\to+\infty}\exp\left\{-\int_{t}^{t+\phi\left(\theta/R\left(y/L(\bar{F}(t)^{-1})\right)\right)}L(R(z))dz\right\}\nonumber\\
&=&\lim_{t\to+\infty}\exp\left\{-\int_{t}^{t+y/L(R(t))}L(R(z))dz\right\}\nonumber\\
&=&e^{-y},\nonumber
\end{eqnarray}
and consequently,
$$\lim_{t\to+\infty}P_{x}\left(\exp\left\{-\theta\,\frac{X_{t}}{R\left(y/L(\bar{F}(t)^{-1})\right)}\right\}|\tau>t\right)=1-e^{-y}.$$
Note that $1-e^{-y}$ is the Laplace transform of
the defective law which assigns mass $1-e^{-y}$ at $0$ and no mass in $(0,+\infty)$.
It follows from the continuity theory
for Laplace transform (see, for example \cite[Section 6.6]{Chung})
that
$$\lim_{t\to+\infty}P_{x}\left(X_{t}\le R\left(y/L(\bar{F}(t)^{-1})\right) |\tau>t\right)=1-e^{-y},$$
or equivalently by \eqref{them4.0}
$$\lim_{t\to+\infty}P_{x}\left(L(\bar{F}(t)^{-1})V(X_{t})\le y
|\tau>t\right)=1-e^{-y}.$$\qed

\section{Examples}

In this section we collect a few examples of branching mechanisms that satisfy the assumptions in Section 3 or Section 4. Branching mechanisms in Examples 1, 2 and 4 are well-known. It follows from \cite[Proposition 5.2]{Ren} that $\psi(\lambda) = \lambda f(\lambda)$ is a critical branching mechanism if and only if $f$ is a Bernstein function
and there exists $b\ge 0$ such that $f(\lambda)=b\lambda+\int^\infty_0(1-e^{-x\lambda})g(x)dx$ with $g\ge 0$ decreasing and
$\int^\infty_0(x\wedge 1)g(x)dx<\infty$.
 Branching mechanisms in Examples 3 and 5 are in given in this from. We refer the reader to \cite{Ren} for more information on the connections between branching mechanisms and Bernstein functions, and \cite{SSV} for more examples of Bernstein functions.

\bigskip
\textbf{Example 1.} Let $\psi(\lambda)=c\lambda^{1+\alpha}$ where
$c>0$ and $\alpha\in (0,1]$. In this case
$\phi(t)=(c\alpha)^{-1}\lambda^{-\alpha}$, $\varphi(t)=(c\alpha
t)^{-1/\alpha}$. Thus we have
$$\bar{F}(t)=1-\exp\{-(c\alpha t)^{-1/\alpha}\}\sim (c\alpha t)^{-1/\alpha}\quad \mbox{as }t\uparrow +\infty.$$
Similarly to \eqref{them1.2}, we get
$$\lim_{t\to+\infty}E_{x}\left(e^{-\theta t^{-1/\alpha}X_{t}}|\tau>t\right)=1-\lim_{t\to+\infty}
\frac{\varphi(t+\phi(\theta
t^{-1/\alpha}))}{\varphi(t)}=1-(1+(c\alpha)^{-1}\theta^{-\alpha})^{-1/\alpha}.$$
Therefore for any $y\ge 0$,
$$\lim_{t\to+\infty}P_{x}\left(t^{-1/\alpha}X_{t}\le y|\tau>t\right)=H_{\alpha}(y),$$
where $H_{\alpha}(y)$ is uniquely determined by its Laplace
transform $$h(\theta)=\int_{0}^{+\infty}e^{-\theta
y}dH_{\alpha}(y)=1-(1+(c\alpha)^{-1}\theta^{-\alpha})^{-1/\alpha}.$$

\begin{remark}
This case was excluded in Pakes et. al. \cite{Pakes2,Pakes3}, and
was studied independently in Haas et.al. \cite{Haas} and Zhang \cite{Zhang}.
More specifically, \cite{Haas} discussed Example 1 as a special case
of self-similar Markov process, while \cite{Zhang} viewed
the corresponding CB process as the scaling limit of a special
sequence of Markov branching processes and exploited limit theorems for some general conditioning events.
\end{remark}

\bigskip
\textbf{Example 2.} If $\psi''(0+)=\sigma<+\infty$, then \eqref{psi}
holds with $\alpha=1$ and $\lim_{s\downarrow 0}L(1/s)=\sigma/2$. By
Karamata's theorem, we have $\phi(z)\sim z^{-1}L(1/z)^{-1}\sim
2/\sigma z$ as $z\downarrow 0$, and $\varphi\in
\mathcal{R}_{-1}(\infty)$. Thus we have
$$\lim_{t\to+\infty}E_{x}\left(e^{-\theta X_{t}/t}|\tau>t\right)=1-\lim_{t\to+\infty}
\frac{\varphi((1+\frac{2}{\sigma}\theta^{-1})t)}{\varphi(t)}=1-(1+\frac{2}{\sigma}\,\theta^{-1})^{-1}.$$
Therefore
$$\bar{F}(t)\sim\frac{2}{\sigma t}\mbox{\quad as }t\uparrow +\infty,$$
and for any $y\ge 0$,
$$\lim_{t\to+\infty}P_{x}\left(X_{t}/t>y|\tau>t\right)=e^{-\frac{2}{\sigma}\,y}.$$
This conditional convergence was proved independently in
Li \cite{Li} and Lambert \cite{Lambert}.

\bigskip
\textbf{Example 3.} Let $\psi(\lambda)=\lambda(\lambda^{-\alpha}+\lambda^{-\beta})^{-1}$ where $0<\beta<\alpha\le 1$. By \cite{SSV} $(\lambda^{-\alpha}+\lambda^{-\beta})^{-1}$ is a Bernstein function, and then $\psi$ is a branching mechanism. Note that $\psi(\lambda)=\lambda^{1+\alpha}L(1/\lambda)$
with $L(z)=(1+z^{-\alpha+\beta})^{-1}$.
By Karamata's theorem, we
have $g(z):=\phi(1/z)=\int_{0}^{z}\xi^{\alpha-1}/L(\xi)d\xi\in
\mathcal{R}_{\alpha}(\infty)$, and
$$g(z)\sim \alpha^{-1}z^{\alpha}L(z)^{-1}\sim \alpha^{-1}z^{\alpha}=:h(z)\mbox{\quad as }z\uparrow +\infty.$$
Both $g$ and $h$ are strictly increasing on $(0,+\infty)$. Let
$g^{-1}$ and $h^{-1}$ respectively denote the inverse functions of $g$ and $h$.
Since
$$1=g(g^{-1}(z))/h(h^{-1}(z))\sim g(g^{-1}(z))/g(h^{-1}(z)),$$ by Lemma \ref{lem1} we have
$g^{-1}(z)\sim h^{-1}(z)=(\alpha z)^{1/\alpha}$ as $z\uparrow +\infty$.
Consequently, $\varphi(t)=1/g^{-1}(t)\sim (\alpha
t)^{-1/\alpha}$ as $t\uparrow +\infty$. Therefore, we have
$$\bar{F}(t)\sim (\alpha t)^{-1/\alpha}\quad\mbox{as }t\to+\infty,$$
and for any $y\ge 0$,
$$\lim_{t\to+\infty}P_{x}\left(t^{-1/\alpha}X_{t}\le y|\tau>t\right)=H_{\alpha}(y),$$
where $H_{\alpha}(y)$ has the Laplace transform
$$h_{\alpha}(\theta)=1-(1+\alpha^{-1}\theta^{-\alpha})^{-1/\alpha}.$$

\bigskip
\textbf{Example 4.} Let
$\psi(\lambda)=\lambda^{1+\beta}+\lambda^{1+\gamma}$,
$0<\gamma<\beta\le 1$. Then
$\psi(\lambda)=\lambda^{1+\gamma}L(1/\lambda)$ with
$L(z)=1+z^{\gamma-\beta}\in\mathcal{S}$.
Using similar arguments as that in Example 3, we have
$$\bar{F}(t)\sim (\gamma t)^{-1/\gamma}\quad\mbox{as }t\to+\infty,$$
and for any $y\ge 0,$
$$\lim_{t\to+\infty}P_{x}\left(t^{-1/\gamma}X_{t}\le y|\tau>t\right)=H_{\gamma}(y),$$
where $H_{\gamma}(y)$ has the Laplace transform:
$$h_{\gamma}(\theta)=1-(1+\gamma^{-1}\theta^{-\gamma})^{-1/\gamma}.$$

\bigskip
\textbf{Example 5.} Let
$\psi(\lambda)=\lambda\log^{-\beta}(1+\lambda^{-1})$, $\beta\in
(0,1]$ and where $\log^{-\beta}(1+\lambda^{-1})$ is a Bernstein function (see \cite[P.133]{Ren}).
Then $\psi$ satisfies \eqref{psi} with $\alpha=0$ and
$L(z)=\log^{-\beta}(1+z)$. Immediately we have $V(z)\sim
(\beta+1)^{-1}\log^{\beta+1}z$ and $L(z)\sim \log^{-\beta}z$ as
$z\uparrow +\infty$.
Inserting the asymptotic equivalents of $V$ and
$L$ into Theorem \ref{them4}, we get
$$-\log\bar{F}(t)\sim [(\beta+1)t]^{\frac{1}{\beta+1}},\mbox{\quad as }t\uparrow +\infty,$$
and
$$\lim_{t\to+\infty}P_{x}\left(\frac{\log^{\beta+1}X_{t}}{(\beta+1)\log^{\beta}(\bar{F}(t)^{-1})}\le y\,|\,\tau>t\right)=1-e^{-y}$$
for any $x>0$ and $y\ge 0$.

\bigskip
\textbf{Acknowledgement} We would like to thank Professor Mei Zhang
and her student Xin Zhang from Beijing Normal University for sending
us their related work on this topic.

\vskip 0.3truein \vskip 0.3truein

 \noindent {\bf Yan-Xia Ren:} LMAM School of Mathematical Sciences $\&$ Center for Statistical Science, Peking
 University, Beijing 100871, P. R. China,\\
 E-mail: {\tt yxren@math.pku.edu.cn} \\

\bigskip
\noindent {\bf Ting Yang:} Institute of Applied Mathematics, Academy
of Mathematics and Systems Science, Chinese Academy of Sciences,
Beijing 100080, P.R.China,\\
 E-mail: {\tt yangt@amss.ac.cn } \\
 \bigskip

 \noindent {\bf Guo-Huan Zhao:} LMAM School of Mathematical Sciences, Peking
 University, Beijing 100871, P. R. China,\\
 E-mail: {\tt 1101110040@math.pku.edu.cn } \\

\end{document}